\documentclass[11pt]{article}
\usepackage{amssymb}
\usepackage{amsmath}
\def\.{\partial_t}
\def\:{\partial_{tt}}
\def\d{\partial}

\def\D{\Delta}

\def\a{\alpha}
\def\e{\varepsilon}
\def\eps{\varepsilon}
\def\f{\varphi}
\def\th{\theta}
\def\l{\lambda}

\def\r{\varrho}
\def\s{\sigma}
\def\t{\tau}

\def\O{\Omega}
\def\G{\Gamma}
\def\ga{\gamma}

\def\Ac{\mathfrak{A}}
\def\cA{\mathcal{A}}
\def\cD{\mathcal{D}}
\def\sD{\mathcal{D}}
\def\cK{\mathcal{K}}
\def\cL{\mathcal{L}}
\def\cE{\mathcal{E}}
\def\H20{\mathcal{H}_{2,0}}
\def\R{\mathbb{R}}

\def\cN{\mathcal{N}}
\def\cH{\mathcal{H}}

\def\<{\langle}
\def\>{\rangle}
\def\8{\infty}
\def\I{\int\limits}

\oddsidemargin
0.1cm \textwidth 16.2cm

\newtheorem{lemma}{Lemma}[section]
\newtheorem{theorem}[lemma]{Theorem}
\newtheorem{remark}[lemma]{Remark}
\newtheorem{proposition}[lemma]{Proposition}

\newtheorem{assumption}[lemma]{Assumption}

\newenvironment{declaration}[1]{\trivlist
\item[\hskip \labelsep{\bf #1 }]\ignorespaces}{\endtrivlist}
\newenvironment{proofof}[1]{\begin{declaration}{#1}}{\hfill
$\square$\end{declaration}}
\newenvironment{proof}{\begin{proofof}{Proof.}}{\end{proofof}}

\begin{document}

\title{Long-time dynamics in plate models \\ with strong nonlinear  damping}
\author{Igor
Chueshov\thanks{e-mail:
chueshov@univer.kharkov.ua}  ~~ and ~~
 Stanislav Kolbasin\thanks{e-mail: staskolbasin@gmail.com}
\\ \\
Department of Mechanics and Mathematics, \\ Kharkov
National University, \\ Kharkov, 61077,  Ukraine
 }
 \maketitle

\begin{abstract}
We study long-time dynamics of a class of
 abstract   second order in time evolution equations
in a Hilbert space with the damping term depending both on displacement
and velocity. This damping    represents the nonlinear strong dissipation
phenomenon perturbed with relatively compact terms.
Our main result states
the existence of a compact finite dimensional attractor.
We study properties of this attractor.
We also establish  the existence of a fractal exponential
 attractor and give the conditions that guarantee the existence
of a finite number of determining functionals.
In the case when the set of equilibria is finite and hyperbolic
we show that every trajectory is attracted by some equilibrium
with exponential rate.
Our arguments   involve a recently
developed   method based on the ``compensated'' compactness and
quasi-stability estimates. As an application we consider the
nonlinear Kirchhoff, Karman and Berger plate models with different
types of boundary conditions and  strong damping terms. Our results
can be also applied to the nonlinear wave equations.
\end{abstract}
\section{Introduction}

We study a class of plate models with the strong nonlinear damping
which abstract form is the following Cauchy problem
in a separable Hilbert space $H$:
\begin{equation}\label{abs-1}
\:u + D(u,\.u)+\cA u+F(u)=0,~~t>0;\quad u|_{t=0}=u_0,  ~~ \.u|_{t=0}=u_1.
\end{equation}
We impose the following set of  hypotheses:
\begin{assumption}\label{A1}
{\rm
\begin{enumerate}
  \item[{\bf (A)}] The operator  $\cA$ is a  linear
self-adjoint positive operator densely
defined on a separable Hil\-bert space $H$  (we denote by $|\cdot |$
and $( \cdot, \cdot )$ the norm and 
the scalar product in this space). We  assume that
the resolvent of $\cA$ is compact in $H$. We also denote
by $H_s$ (with $s>0$) the domain $\sD(\cA^{s/2})$ equipped with the graph norm
$|\cdot |_s=|\cA^{s/2}\cdot|$. In this case $H_{-s}$ denotes the completion
of $H$ with respect to the norm $|\cdot |_{-s}=|\cA^{-s/2}\cdot|$.
Below we  denote by $\{e_k\}$ the orthonormal basis in $H$ consisting of the eigenfunctions
of the operator $\cA$ : 
\[
\cA e_k=\l_k e_k,
\quad 0<\l_1\le \l_2\le\cdots, \quad \lim_{k\to\infty}\l_k =\8,
\]
and by    $P_N$  the orthoprojector
onto Span${}\{ e_k:\; k=1,2,\ldots, N \}$.

  \item[{\bf (D)}] For some value $ \th \in(0,1]$ the damping
operator $D$ maps  $H_1\times H_\th$ into $H_{-\th}$ and
possesses the properties:
\begin{enumerate}
    \item[(i)] For every $\varrho>0$ there exist $\a_\r>0$ and $\beta_\r>0$
such that
\begin{equation}\label{K-damp1}
(D(u, v),v)\ge  \a_\varrho  | v|_\th^2 ~~\mbox{and}~~ |D(u, v)|_{-\th}\le  \beta_\varrho | v|_\th
\end{equation}
for all $(u;v)\in H_1\times H_{\th}$, $|u |_1^2 +|v|^2 \leq \varrho^2$
(it is allowed that $\a_\r\to 0$ and $\beta_\r\to\infty$
as $\r \rightarrow \8$).
  \item[(ii)] For every $\varrho>0$ there exist $\ga_\r>0$ and $C_\r>0$
such that
\begin{equation}\label{K-damp1a}
(D(u^1, v^1)-D(u^2, v^2) ,v^1-v^2)\ge  \ga_\varrho  | v^1-v^2|_\th^2 -
C_\r |u^1-u^2|^2_{1-\delta} (1+ | v^1|^2_\th + | v^2|^2_\th )
\end{equation}
for all $(u^i;v^i)\in H_1\times H_{\th}$, 
$|u^i |_1^2 +|v^i|^2 \leq \varrho^2$,
and for some $\delta>0$
(it is allowed that $\ga_\r\to 0$ and $C_\r\to\infty$
as $\r \rightarrow \8$).
 \item[(iii)] We assume that for every $u\in L_\infty(0,T;H_1)$
the mapping  $v\mapsto D(u, v)$
is weakly continuous from $L_2(0,T;H_\th)$ into  $L_2(0,T;H_{-\th})$
and
\[
 | D(u^1,v )-D(u^2,v)|_{-1} \le  C_\varrho
  |u^1-u^2|_{1-\delta}  (1+ | v|^2_\th)
\]
for all $|u^i |_1, |v|  \leq \varrho$, $u^i\in H_1$, $v\in H_\th$, and
for  some $\delta>0$.
 \item[(iv)] In the case $\th<1$ we  have
\begin{equation}\label{K-damp4-a}
 | D(u^1,v^1 )-D(u^2,v^2)|_{-\th} \le  C_\varrho \left[ |v^1-v^2|_\th+
  |u^1-u^2|_{1}  (1+ | v^1|_\th+| v^2|_\th)\right]
\end{equation}
for all $(u^i;v^i)\in H_1\times H_{\th}$,
$|u^i |_1^2 +|v^i|^2 \leq \varrho^2$.
In the case $\th=1$ we assume a stronger
(compared to (\ref{K-damp4-a}) with $\th=1$) inequality:
\begin{equation}\label{K-damp4-b}
 | D(u^1,v^1 )-D(u^2,v^2)|_{-1} \le  C_\varrho \left[|v^1-v^2|_1+
  |u^1-u^2|_{1-\delta}  (1+ | v^1|_1+| v^2|_1)\right]
\end{equation}
for all $(u^i;v^i)\in H_1\times H_{1}$,
$|u^i |_1^2 +|v^i|^2 \leq \varrho^2$, where $\delta>0$.
\end{enumerate}
We note that the conditions in (\ref{K-damp1}) follows from
(\ref{K-damp1a}) and (\ref{K-damp4-a}) or (\ref{K-damp4-b})
provided  $D(u,0)\equiv 0$ for every $u\in H_1$.
  \item[{\bf (F)}] There exists $\delta>0$ such that
the nonlinear operator $F$  maps $H_{1-\delta}$ into $H_{-\th}$ and is
locally Lipschitz, i.e.,
\begin{equation}\label{8.1.1c}
|F(u_1) -F(u_2)|_{-\th}~\leq~L(\varrho) | u_1-u_2|_{1-\delta},\quad
\forall |\cA^{1/2} u_i |  \leq \varrho\,,
\end{equation}
We also assume that $ F(u)= \Pi^\prime(u)$, where $\Pi(u)$ is a
$C^1$-functional on $H_1=\sD (\cA^{1/2})$, and 
${}^\prime$ stands for the Fr\'echet derivative.
We assume that $\Pi(u)$
 is  locally bounded on $H_1$ and there exist $\eta<1/2$
 and $C\ge 0$ such that
\begin{equation}\label{8.1.1c1}
 \eta |\cA^{1/2}u|^2 +\Pi(u)+C \ge 0\;,\quad u\in H_1=\sD ( \cA^{1/2})\;.
\end{equation}
\end{enumerate}
}
\end{assumption}
\begin{remark}
{\rm Our conditions can be relaxed in different directions. For instance, for
the well-posedness in Theorem~\ref{pr:wp1} we can assume that $\delta=0$
in (\ref{K-damp1a}) and (\ref{8.1.1c}), and  require $F$ to
map $H_{1-\s}$ into $H_{-l}$ continuously for some $\s>0$ and $l>0$. In a similar way,
 instead of (Diii) we can assume other continuity properties of the damping
operator $D$ that will allow us to perform the limit transition in
the corresponding Galerkin approximations.
To obtain the global attractor existence, we can also relax (\ref{K-damp1a}) by including into its
latter term   the expression $\eps |u^1-u^2|_1$ with a small parameter  $\eps$
(see, e.g., Assumption~3.21 in \cite{ChuLas} for a similar requirement in
the case of the monotone damping).
Moreover, instead of
(\ref{K-damp1a}) we can assume that
\begin{equation}\label{new-req}
(D(u^1, v^1)-D(u^2, v^2) ,v^1-v^2)\ge  \ga_\varrho  | v^1-v^2|_\th^2 -
|u^1-u^2|^2_{1} \left[\eps + C_\r(\eps) \left(| v^1|^2_\th + | v^2|^2_\th
\right)\right]
\end{equation}
for every $\eps>0$.
 However, we do not pursue
these possible generalizations because our abstract hypotheses are
motivated by the plate models described below.
We note that in this paper we concentrate on the case
when $\theta$ is positive.
For some results for the case when $\th=0$ and
$D(u,u_t)$ is linear with respect to $u_t$, we refer to \cite{ChuKol}.
It is also worth  mentioning that our damping operator $D$ is positive (see
the first relation in  (\ref{K-damp1})) but not monotone
(see (\ref{K-damp1a})) in general. Thus, we cannot apply here the theory
 developed in \cite{ChuLas}.
}
\end{remark}

Our main applications are plate models (with hinged boundary conditions,
for definiteness).
In this case
\begin{itemize}
    \item
$\cA=(-\Delta_D)^2$,  where
$\Delta_D$ is the Laplace operator in a bounded smooth domain
$\O$ in $\R^2$ with the Dirichlet boundary conditions.
We have then that  $H=L_2(\O)$ and
$$ \sD ( \cA)=\left\{ u\in H^4(\O) \, : \;
 u=\Delta u =0 ~~{\rm on}~~\d\O\right\}.
$$
We also have that $H_{s}=(H^{2s}\cap H^1_0)(\O)$
for $1/2\le s\le 1$ and  $H_{s}=H^{2s}_0(\O)$ for $0<s<1/2$.
Here above $H^\s(\O)$ is the corresponding $L_2$-based Sobolev space,
$H^\s_0(\O)$ is the completion of $C^\infty_0(\O)$ in $H^\s(\O)$.
\item  The damping operator $D(u,u_t)$ may have the form
\begin{equation}\label{damping}
D(u,u_t)=\Delta\left[\s_0(u)\Delta u_t\right]-
{\rm div}\,\left[\s_1(u,\nabla u)\nabla u_t\right]+g(u,u_t),
\end{equation}
where $\s_0(s_1)$, $\s_1(s_1,s_2,s_3)$ and $g(s_1,s_2)$
are locally Lipschitz functions of $s_i\in\R$, $i=1,2,3$, such that $\s_0(s_1)>0$, $\s_1(s_1,s_2,s_3)\ge 0$ and $g(s_1,s_2)s_2\ge 0$.
Also the functions  $\s_1$ and $g$
satisfy  some growth conditions (for a more detailed discussion
of properties of the damping functions we refer to Section~\ref{s:appl}
below). We note that every term in (\ref{damping})
represents a different type of damping mechanisms. The first one is the so-called viscoelastic Kelvin--Voight damping, the second one represents
the structural damping and the term  $g(u,u_t)$ is the dynamical friction
(or viscous damping). We refer to \cite[Chapter 3]{lt-book}
and to the references therein for a discussion of stability
properties caused by each type  of the damping terms   in the case of
linear systems.
\item The nonlinear feedback (elastic) force $F(u)$ may have one
of the following forms (which represent different plate models):
\begin{itemize}
  \item[(a)] {\sl Kirchhoff model}: $F(u)$ is the Nemytskii operator
\begin{equation}
    \label{kirch-force}
u\mapsto - \kappa\cdot {\rm div}\left\{|\nabla u|^q\nabla u
-\mu |\nabla u|^r\nabla u
\right\}+ \f(u)-p(x),
\end{equation} 
 where  $\kappa\ge 0$, $q>r\ge 0$, $\mu\in \R$ are parameters, $p\in L_2(\O)$ and
$\f\in {\rm Lip_{loc}}(\R)$ fulfills the  condition
 \begin{equation}
    \label{phi_condition}
    \underset{|s|\to\8}{\liminf}{\frac{\f(s)}{s}}>-\l_1^2,
\end{equation} 
 where $\l_1$ is the first eigenvalue of the
Laplacian with the Dirichlet boundary conditions.
  \item[(b)] {\sl Von Karman model:} $F(u)=-[u, v(u)+F_0]-p(x)$, where
 $F_0\in H^4(\O) $ and $p\in L_2(\O)$ are given functions,
the von Karman bracket $[u,v]$  is given by
\begin{equation*}
[u,v] = \partial ^{2}_{x_{1}} u\cdot \partial ^{2}_{x_{2}} v +
\partial ^{2}_{x_{2}} u\cdot \partial ^{2}_{x_{1}} v -
2\cdot \partial ^{2}_{x_{1}x_{2}} u\cdot \partial ^{2}_{x_{1}x_{2}}
v ,
\end{equation*} and
the Airy stress function $v(u) $ solves the following  elliptic
problem
\begin{equation}\label{airy-1}
\Delta^2 v(u)+[u,u] =0 ~~{\rm in}~~  \Omega,\quad \frac{\d v(u)}{\d
n} = v(u) =0 ~~{\rm on}~~  \d\O.
\end{equation}
Von  Karman equations are well known in nonlinear elasticity and
constitute a basic model describing nonlinear oscillations of a
plate accounting for  large deflections, see \cite{Lio69,cl-book}  and the
references therein.
  \item[(c)] {\sl Berger Model:} In this case the feedback force has the form
  $$
  F(u)=- \left[ \kappa \int_\O |\nabla u|^2 dx-\G\right] \D u -p(x),
$$
where $\kappa>0$ and $\G\in\R$ are parameters,  $p\in L_2(\O)$;
 for some details and  references see, e.g.,
\cite[Chapter 4]{Chueshov} and \cite[Chapter 7]{ChuLas}.
\end{itemize}
\end{itemize}

Long-time dynamics of second order equations with a nonlinear damping
was studied by many authors. We refer to \cite{LasieckaBarbuRammaha, GattiPata_1D, PataZelik_2D, PataZelik_3D, ChuKol,Kolbasin} for the
case of a damping   with a
nonlinear displacement-dependent coefficient   and
to  \cite{ChuLas_JDDE_2004,ChuLas,cl-book}
and to the references therein
for the velocity-dependent damping.
Models with different types of a strong (linear) damping
in wave equations were considered in
\cite{CaCh02,ChDl06,KaZe09,PataZelik06,YaSu09},
see also the literature quoted in these references.
\par

 The main novelty of the current paper is the
following: (i) we can consider the strong nonlinear displacement- and velocity-dependent damping of a general structure
(thus we cannot use analyticity of the corresponding
model with a zero source term which takes place in the case when $D(u,u_t)=Bu_t$,
where $B$ is a self-adjoint operator satisfying (\ref{K-damp1}) with
$\th\in [1/2,1]$, see, e.g.,  \cite[Chapter 3]{lt-book} and the
references therein);
(ii) this damping can be perturbed by low order terms.
\par
Our main result (see Theorem~\ref{th:attractor}) states the existence of a
compact global attractor and describes other asymptotic (long-term)
properties of the system generated by (\ref{abs-1}). To
establish this result we use the recently developed approach (see \cite{ChuLas_JDDE_2004} and also  \cite{ChuLas} and
    \cite[Chapters 7,8]{cl-book}). We first prove that the corresponding system is quasi-stable in the sense of the  definition given in
\cite[Section 7.9]{cl-book},
and then we apply general theorems on properties of quasi-stable
systems. In the same framework 
we also establish a result on the rate of 
stabilization (see Theorem~\ref{th:rate}) which 
states that under some additional conditions every solution 
is attracted by an equilibrium with an exponential rate. 
To obtain this result  we rely on
some type of  the observability inequality and
 use the same idea as in \cite[Section 4.3]{ChuLas} (see also \cite{ChuLas_JDDE_2004} and \cite{cl-book}).

\par
The paper is organized as follows. In the preliminary Section~\ref{sect2}
we discuss well-posedness of our abstract model and the dynamical system generation. We also recall several notions and results from the theory
of dissipative dynamical systems. Our main results on 
the global attractor for (\ref{abs-1})
and on asymptotic properties of individual trajectories
are stated
in Section~\ref{s:main-r}. The proofs  based on
the quasi-stability property of the corresponding system are given in
Section~\ref{sect4}. In Section~\ref{s:appl} we discuss some applications.
\par
Below  constants denoted by the same symbol may vary from line to line.

\section{Preliminaries}\label{sect2}
In this section we show that problem (\ref{abs-1}) generates a
dynamical system.

\subsection{Well-posedness}\label{sect3}

We first prove the  existence and uniqueness of weak solutions
 to problem (\ref{abs-1}). We  recall that a function $u(t)$
is  a weak solution
to (\ref{abs-1}) on an interval $[0,T]$ if
\[
u \in L_\8(0,T; \sD(\cA^{1/2})), \quad \.u \in L_\8(0,T; H)\cap L_2(0,T;\sD(\cA^{\th/2}))
\]
and  (\ref{abs-1}) is satisfied in the sense of distributions.
\par
The main statement of the section is the
following assertion which also contains some
auxiliary solution properties
needed for the results on  asymptotic dynamics.

\begin{theorem}\label{pr:wp1}
Let Assumption~\ref{A1} be in force and
 $(u_0; u_1) \in \cH\equiv  \cD(\cA^{1/2})\times H$.
Then the following assertions hold.
\begin{enumerate}
  \item[{\bf 1.}]
  Problem (\ref{abs-1}) has a unique weak solution
$u(t)$ on $\R_+$. This
solution belongs to the class
\[
{\cal W}\equiv C(\R_+;\sD ( \cA^{1/2}))\cap C^1(\R_+;  H ),
\]
and the following energy relation
\begin{equation}\label{8.1.4}
\cE(u(t), \.u(t))+\int_0^t (D(u(\tau), \d_\t u(\tau)),\d_\t u(\tau)) d\tau=
\cE(u_0, u_1)
\end{equation}
holds for every $t>0$, where the energy $\cE$ is defined by the formula
\begin{equation*}
\cE(u_0, u_1)=E(u_0, u_1)+\Pi(u_0)
\equiv\frac12\left( |u_1|^2 + \big|\cA^{1/2} u_0\big|^2\right)+\Pi(u_0).
\end{equation*}
Moreover, this solution $u(t)$  satisfies the estimate
\begin{equation}
\label{.u-int-bound}
\sup_{t\ge 0}E(u(t), \.u(t))+\I_0^{+\8}|\cA^{\th/2}\.u(t)|^2dt \leq C(R)
~~if~~ E(u_0, u_1)\le R^2.
\end{equation}

  \item[{\bf 2.}]
If $u^1(t)$ and $u^2(t)$ are two weak solutions
such that $E(u^i(0), \.u^i(0))\le R^2$, $i=1,2$, then
their difference $z(t)=u^1(t)-u^2(t)$ satisfies the relation
\begin{equation}
\label{dif-bnd}
E(z(t), \.z(t))+a_R\I_0^{t}|\cA^{\th/2}\.z(\t)|^2d\t \leq b_R
E(z(0), \.z(0))e^{c_Rt}
\end{equation}
for some constants $a_R,b_R, c_R>0$.

\end{enumerate}

\end{theorem}
\begin{proof}
To prove the existence of solutions, we use the standard Galerkin method of seeking for approximations of the form
$$
u_N(t)=\sum\limits_{k=1}^N C_k(t)e_k,\quad N=1,2,\ldots
$$
that solve the finite-dimensional projections of (\ref{abs-1}). Such solutions exist, and after multiplication of the corresponding projection of
 (\ref{abs-1}) by $\.u_N(t)$ we get that $u_N(t)$ satisfies the energy relation \eqref{8.1.4}. By \eqref{8.1.1c1} we obtain that
\begin{equation*}
c_0 E(u_0, u_1)-c_1\le \cE(u_0, u_1)\le C(R)
\end{equation*}
whenever $E(u_0, u_1)\le R^2$.
Therefore, by (\ref{K-damp1})  the energy relation for  $u_N(t)$ yields
 estimate \eqref{.u-int-bound} for approximate solutions
with the constant $C(R)$ independent of $N$.
Using the equation for $u_N(t)$ and also the conditions (\ref{K-damp1}) and (\ref{8.1.1c}),
 it can be shown in the standard way that
\[
\I_0^T |\cA^{-1/2}\:u_N(t)|^2 dt \le C_T(R), \quad N=1,2,\ldots,
\]
for every $T>0$.
These a priori estimates show that
$(u_N;\.u_N;\:u_N)$ is ${}^*$-weakly compact in
\[
L_\infty(0,T; H_1)\times\left[ L_\infty(0,T; H)\cap L_2(0,T; H_\th)\right]
\times L_2(0,T; H_{-1})\quad \mbox{for every}~~ T>0.
\]
Thus the Aubin-Dubinsky theorem (see \cite[Corollary 4]{Simon})
yields that $(u_N;\.u_N)$ is  compact in
$C(0,T; H_{1-\e}\times H_{-\e})$
for every $\e>0$.
These compactness properties make it possible to show the existence
of weak solutions satisfying \eqref{.u-int-bound}. For the limit transition in the nonlinear terms we use the property (Diii) in Assumption~\ref{A1}
and relation \eqref{8.1.1c}. It is also clear that $t\mapsto (u(t);\.u(t))$
is a weakly continuous function in $\cH=\sD(\cA^{1/2})\times H$.
To obtain the energy relation in (\ref{8.1.4}), we note that the function
$u^n(t)=P_nu(t)$  solves an equation of the form
\[
\:u^n+\cA u^n=h(t)
\]
with some $h\in L_2(0,T;H)$. This makes it possible to obtain a certain energy relation for $u^n$ which gives us  (\ref{8.1.4}) after the limit transition
$n\to\infty$. This also allows us to obtain the strong continuity properties
of $t\mapsto (u(t);\.u(t))$ in $\cH$ by the standard method.
\par
To prove \eqref{dif-bnd}, we note that $z(t)=u^1(t)-u^2(t)$
solves the equation
\begin{equation}\label{abs-dif}
    \:z + D(u^1,\.u^1)- D(u^2,\.u^2) +\cA z  +F(u^1)-F(u^2)=0.
\end{equation}
Thus, multiplying this equation by $\d_t z$ and integrating
from $s$ to $t$ we have
\begin{equation}
\label{abs-dif-z_t-0}
    E_z(t)+\I_s^t (D(u^1,\d_\t u^1)- D(u^2,\d_\t u^2),\d_\t z)d\t  = E_z(s)
    - \I_s^t \big( F(u^1)-F(u^2),\d_\t z \big)d\t
\end{equation}
for any $0\le s<t$, where $E_z(t)=E(z(t), \.z(t))$. Therefore, using \eqref{K-damp1a},
\eqref{8.1.1c}  and
\eqref{.u-int-bound} we obtain that
\[
    E_z(t)+\frac{\ga_R}2\I_s^t \big|\cA^{\th/2}\d_\t z\big|^2 d\t  \le  E_z(s)
    + c_R \I_s^t \big( 1+ \big|\d_\t u_1\big|_\th^2 + \big|\d_\t u_2\big|_\th^2  \big) | \cA^{1/2}z|^2d\t, ~~ s<t,
\]
for some $c_R>0$.
Now we can apply the Gronwall lemma to obtain \eqref{dif-bnd} which, in particular, implies
the uniqueness of weak solutions.
\end{proof}

\subsection{Generation of a dynamical system}
We recall that a {\em dynamical system} (see, e.g., \cite{Chueshov,Hale,Temam}) is a pair
$\big(X, S(t)\big)$ of a complete metric space $X$ and a family of
continuous mappings $S(t):X\mapsto X,\ t\ge 0$,
such that (i) $t\mapsto S(t)y$ is continuous in $X$ for every $y\in X$
and (ii) the semigroup property is satisfied, i.e.,
 $ S(t+\t) = S(t)\circ S(\t)$ for any $t,\t\ge 0$ and
$S(0)$ is the identity operator.
\par
We also recall that the system $\big(X, S(t)\big)$   is gradient
if it  possesses a
\textit{strict Lyapunov function}, i.e.,   there exists a
continuous functional $\Phi(y)$ on $X$ such that
(i) $\Phi\big(S(t)y\big) \leq \Phi(y)$  for all $t\geq 0$ and $y\in X$;
(ii) the equality $\Phi(y)=\Phi(S(t)y)$ may take place for all $t>0$ if only $y$
is a stationary point of $S(t)$. 
\par
Applying Theorem \ref{pr:wp1} we obtain the following assertion.
\begin{proposition}\label{pr:generation}
Let Assumption~\ref{A1} be in force. Then
 problem \eqref{abs-1} generates a dynamical system in
the space $\cH=\sD(\cA^{1/2})\times H$ with the evolution operator $S(t)$
given by
\[
S(t)y=(u(t);\.u(t)),\quad
\mbox{where $y=(u_0;u_1)$ and $u(t)$ solves (\ref{abs-1}).}
\]
This system is gradient with  the full energy $\cE(u_0;u_1)$ as a strict Lyapunov function (this follows from  the energy relation in \eqref{8.1.4}).
\end{proposition}
We also recall  that a
 system $\big( X, S(t)\big)$ is called \textit{asymptotically smooth}
(see \cite{Hale})
if for any  closed bounded set $B\subset X$ that is positively invariant ($S(t)B\subseteq B$)
one can find a compact set $\cK=\cK(B)$ which uniformly attracts $B$:
$\sup\{{\rm dist}_X(S(t)y,\cK):\ y\in B\}\to 0$ as $t\to\infty$.
The \textit{global attractor} (see, e.g., \cite{BabinVishik, Chueshov,Hale,Temam}) of a dynamical system  $\big( X, S(t)\big)$
is defined as a bounded closed  set $\Ac\subset X$
which is  invariant ($S(t)\Ac=\Ac$ for all $t>0$) and  uniformly  attracts
all other bounded  sets:
$$
\lim_{t\to\8} \sup\{{\rm dist}_X(S(t)y,\Ac):\ y\in B\} = 0
\quad\mbox{for any bounded  set $B$ in $X$.}
$$
In this paper we use the following criterion of the global attractor existence for gradient systems (see, e.g., \cite[Theorem 4.6]{Raugel}):
\begin{theorem}\label{Theorem 2.2.}
Let $\big( X, S(t) \big)$ be an asymptotically smooth gradient system
such that for any bounded set
 $B\subset X$ there exists $\t>0$ such that
 $\ga_\t(B)\equiv\bigcup_{t\geq\t}S(t)B$ is bounded. If the set $\cN$ of
stationary points  is bounded, then $\big( X, S(t) \big)$ has a
compact global attractor $\Ac$ which coincides with the unstable set $\mathbb{M}_+(\cN)$
emanating from $\cN$, i.e., $\Ac=\mathbb{M}_+(\cN)$.
\end{theorem}
We recall (see, e.g., \cite{BabinVishik})
 that the \textit{unstable set} $\mathbb{M}_+(\cN)$ emanating from $\cN$ is a subset
 of $X$ such that for each $z\in\mathbb{M}_+(\cN)$
there exists a full trajectory $\{y(t): t\in\R\}$ satisfying
$u(0) = z$ and ${\rm dist}_X(y(t),\cN) \to  0$ as $t\to -\8$.

\begin{remark}\label{re:th2.2}
{\rm
We note that we can avoid the hypothesis of Theorem~\ref{Theorem 2.2.} that $\ga_\t(B)$ is bounded
if  the following requirements on the corresponding Lyapunov function $\Phi(x)$ are added: 
(i) $\Phi(x)$ is bounded from above on any bounded set; 
(ii) the set $\Phi_R=\{x\in X: \Phi(x)\le R\}$ is bounded for
every $R$ (see, e.g., \cite[Corollary 2.29]{ChuLas}).
}
\end{remark}

\section{Main results}\label{s:main-r}
Our first main result is the following theorem.
\begin{theorem}\label{th:attractor}
Let Assumption~\ref{A1}  be in force.
Then the dynamical system $(\cH, S(t))$ 
generated by (\ref{abs-1})
possesses a compact global attractor $\Ac$.
 Moreover,
\begin{enumerate}
    \item[{\bf (1)}] $\Ac=\mathbb{M}_+(\cN)$, where
 $\cN=\{ (u;0)\in\cH : \cA u +F(u)=0\}$ is the set of stationary points and 
\begin{equation}\label{n-stab}
{\rm dist}_\cH (S(t)y,\cN)\equiv
 \inf\left\{ |S(t)y- e|_\cH\, : e\in\cN \right\}\to 0~~as~~ t\to+\infty~~
for~every~y\in\cH.
\end{equation}
\item[{\bf (2)}] This attractor has a finite  fractal dimension.
\item[{\bf (3)}] Any trajectory $\gamma=\{ (u(t); \.u(t)): t\in \R\}$ from the attractor $\Ac$
possesses the property
\begin{equation}\label{u-smth}
(u;\.u;\:u)\in L_\infty (\R;  H_{2-\th}\times H_1\times H),
\end{equation}
and there is $R>0$ such that
\begin{equation}\label{u-smth2}
\sup_{\ga \subset \Ac}\sup_{t\in\R}\left( |u|^2_{2-\th}+|\.u|^2_1+ |\:u|^2\right)\le R^2.
\end{equation}
\item[{\bf (4)}] The system $(\cH, S(t))$ possesses a (generalized)
fractal exponential attractor $\Ac_{exp}$ whose dimension is finite in
the space $\tilde \cH=H_\th\times H_{-1}$.
\item[{\bf (5)}]

 Let   ${\cal L} =  \{ l_j : j= 1,...,N\}$ be a finite set of
functionals on $H_1$ with the completeness
defect
\[
\epsilon_{\cal L}=\epsilon_{\cal L}(H_1,H)\equiv
\sup\left\{ |u|\,:\, u\in H_1, l_j(u)=0,
j=1,...,N,
|u|_1\le1\right\}.
\]
Then  there exists $\eps_0>0$ such that
under the condition $\eps_\cL\le\eps_0$
 the set $\cL$ is (asymptotically)
determining in the sense   that the property
\[
\lim_{t\to\infty}\max_{j} \int_t^{t+1}|l_j(u^1(s)-u^2(s))|^2ds=0
\quad\mbox{for  two solutions $u^1$ and $u^2$}
\]
implies that $\lim_{t\to\infty} |S(t)y_1-S(t)y_2|_\cH=0$.
Here above $S(t)y_i=(u^i(t); \.u^i(t))$, $i=1,2$.
\end{enumerate}
\end{theorem}
We recall that
the {\em fractal dimension} $\dim^X_f M$ of a compact set $M$ in a complete
metric space $X$ is defined as
\[
\dim^X_fM=\limsup_{\eps\to 0}\frac{\ln N(M,\eps)}{\ln (1/\eps)}\;,
\]
where $N(M,\eps)$ is the minimal number of closed sets of
diameter $2\eps$ in $X$ needed to cover the set~$M$.
\par
We also recall (see, e.g., \cite{EFNT94} and
 also \cite{ChuLas,MiZe2008} and the references therein)
that a
 compact set $\Ac_{\rm exp}\subset \cH$
is said to be a (generalized) fractal exponential
attractor
for  the dynamical system $(\cH, S(t))$  iff  $\Ac_{\rm exp}$ is a positively invariant set
of finite fractal dimension (in some extended space $\tilde\cH$) and
for every bounded set $D\subset \cH$ there exist positive constants
$t_D$, $C_D$ and $\gamma_D$ such that
\begin{equation*}
d_\cH\{S(t)D\, |\, \Ac_{\rm exp}\}\equiv
\sup_{x\in D} \mbox{dist}\,_\cH  (S(t)x,\, \Ac_{\rm exp})\le C_D\cdot e^{-\gamma_D(t-t_D)},
\quad t\ge t_D.
\end{equation*}
As for the   determining functionals,
we mention that this notion   goes back
to the papers  by Foias and Prodi~\cite{FP68}
and by Ladyzhenskaya~\cite{Lad75} for the 2D Navier-Stokes equations.
For the further  development of the theory  we refer to
 \cite{cjt97} and to the survey \cite{Chu98}
 and to the references quoted therein (see also \cite[Chap.5]{Chueshov}).
We note that for the first time determining functionals for 
second order (in time) evolution equations with a nonlinear damping was considered in 
\cite{CHuKal2001}, see also a discussion in \cite[Section 8.9]{cl-book} 
We also refer  to \cite{Chu98,Chueshov} for a description of sets
of functionals with a small completeness defect.
Determining modes and nodes are among them.
\medskip\par  
Using the same idea as in \cite{ChuLas_JDDE_2004,ChuLas,cl-book} we can 
establish the following result on convergence of
individual solutions to equilibria  with an exponential rate. 
\begin{theorem}\label{th:rate}
In addition to  Assumption~\ref{A1}  we assume that 
 $F(u)$
is Fr\'echet differentiable and its derivative $F'(u)$ possesses the
properties
\begin{equation}\label{8.6.1a}
|\langle F'(u), w\rangle|_{-1}\le C_R
| w|_{1},~~ w\in H_1,
\end{equation}
and
\begin{equation}\label{8.6.1b}
|\langle F'(u)- F'(v), w\rangle|_{-1}\le C_R
|u-v|_{1-\delta}\cdot |w|_1,~ w\in H_1,
\end{equation}
for any $u,v\in H_1$ such that $|u|_1\le R$ and
$|v|_1\le R$ with $\delta>0$. Here $\langle F'(u), w\rangle$ is the value of
$F'(u)$ on the element $w$.
Let the set $\cN$ be finite and all equilibria be
hyperbolic in the sense that the equation
$\cA u+\langle F'(\phi), u\rangle=0$ has only a trivial solution for
each $(\phi;0)\in\cN$.
Then for any $y\in \cH$ there exists an equilibrium
$e=(\phi;0)\in\cN$ and  constants  $\gamma>0$, $C>0$
such that
\begin{equation}\label{8.6.2}
| S(t)y -e|_\cH\le C e^{-\gamma t},~~ t>0.
\end{equation}
\end{theorem}
We note that this type of  stabilization theorems is well-known in 
literature for different classes of gradient systems, and several approaches 
to the question on  stabilization rates are available 
(see, e.g., \cite{BabinVishik} and also 
 \cite{ChuLas_JDDE_2004,ChuLas,cl-book}
and the references therein). The approach presented 
in \cite{BabinVishik} relies on the analysis of linearized 
dynamics near each equilibrium and requires 
the hyperbolicity  condition  in  a dynamical form.
Here  we use the method developed in
 \cite{ChuLas_JDDE_2004,ChuLas} (see also a discussion in  \cite{cl-book}),
and we need this condition in a weaker form.

\section{Proofs}\label{sect4}
As it was already mentioned, the main ingredient of the proof of Theorem~\ref{th:attractor}
is a quasi-stability  property of the dynamical system $(\cH, S(t))$
generated by (\ref{abs-1}).  
\subsection{Quasi-stability}
We show that under the conditions listed in Assumption~\ref{A1} the  system $(\cH,S(t))$
is quasi-stable in the sense of the definition given
in \cite[Section 7.9]{cl-book}.
Namely, we prove the following proposition.

\begin{proposition}\label{pr:q-st} Let Assumption \ref{A1} hold.
Assume that $u^i(t)$, $i=1,2$ are two weak solutions to problem  (\ref{abs-1}) with
 initial data $y_i=(u^i_0;u^i_1)$ such that $\big|\cA^{1/2}u^i_0\big|^2 +|u^i_1|^2 \leq R^2$ for some $R>0$.
We denote  $S(t)y_i=(u^i(t); \.u^i(t))$, $i=1,2$.
Then there exist $C(R), \ga(R)>0$ such that
\begin{equation} \label{q-stab-est}
   |S(t)y_1-S(t)y_2|^2_\cH\le   C(R)\left[ |y_1-y_2|^2_\cH e^{-\ga(R) t} +
    \int_0^t e^{-\ga(R)(t-\t)}|u^1(\t)-u^2(\t)|^2 d\t \right],\quad t>0.
\end{equation}
\end{proposition}
This type of estimates  was originally introduced in \cite{ChuLas_JDDE_2004}
and related to a decomposition of the evolution operator $S(t)$ into
uniformly exponentially stable and compact parts,
see also a discussion in \cite{ChuLas} and \cite[Section 7.9]{cl-book}.

We start with two preliminary lemmas.

\begin{lemma}\label{le:new-dim1}
Under Assumption~\ref{A1} there exist $T_0>0$ and  a constant
$c>0$ independent of
$T$ such that for any pair $u^1$ and $u^2$ of weak solutions to
(\ref{abs-1}) we have
the following  relation
\begin{eqnarray*}
 TE_z(T) +
\int_{0}^{T} E_z(t) dt
\le c\left\{ \int_0^T |z_t(t)|^2 dt+  \int_0^T |(D(t),z_t)| dt\right.
 \nonumber \\
 + \left. \int_0^T \left|(D(t),z)\right| dt+\Psi_T(u^1,u^2)\right\} 
\end{eqnarray*}
for  every $T\ge T_0$,
where $z(t)=u^1(t)-u^2(t)$, and the functionals $E_z$, $D$ and $\Psi_T$ are defined as 
\begin{eqnarray*}
E_z(t) &= & E_0(z(t),z_t(t))=
\frac12\left( (z_t(t),z_t(t))+(\cA z(t),z(t))\right),
\nonumber
\\
D(t)& = & D(u^1(t),u^1_t(t))- D(u^2(t),u^2_t(t)),
\nonumber \\
\Psi_T(u^1,u^2) &= & \left|  \int_{0}^{T}
(G(\tau), z_t(\tau) ) d\tau\right|
+  \left| \int_{0}^{T}(G(t), z(t) ) dt\right|\nonumber \\
& &  +   \left|  \int_{0}^{T}dt \int_{t}^{T}
(G(\tau), z_t(\tau) ) d\tau \right|
\end{eqnarray*}
with   $G(t)= F ( u^1(t)) - F(u^2(t))$.
\end{lemma}
\begin{proof}
We use the standard arguments involving the multipliers $z_t$ and $z$  for (\ref{abs-dif}). 
We refer to the proof of Lemma 3.23 in \cite{ChuLas}
and also \cite[Lemma 8.3.1]{cl-book}, where this lemma is
proved under another set of hypotheses concerning the damping operator.
However the corresponding argument    does not depend on a structure
of the damping operator. 
\end{proof}

\begin{lemma}\label{le:max}
Let $u^1$ and $u^2$ be  two solutions to (\ref{abs-1})
with the initial data
$(u_0^i;u_1^i)$. We assume that
$|u_1^i|^2 + |u_0^i|_1^2 \leq R^2$. Then
\begin{equation*}
    \underset{[0,T]}{\max}E_z(t) \leq c_0\left[ E_z(T)+
  \int_0^T |(D(t),z_t)| dt\right] + C_R\left[  \int_0^T|\cA^{\th/2}z_t|^2dt +   \int_0^T\big| \cA^{1/2}z \big|^2dt\right],
\end{equation*}
with $z=u^1-u^2$, where $E_z$ and $D(t)$ are the same
as in Lemma~\ref{le:new-dim1}.
\end{lemma}
\begin{proof}
This follows from (\ref{abs-dif-z_t-0}), the Lipschitz property for $F$
and the uniform estimate (\ref{.u-int-bound}) for $u^i(t)$, $i=1,2$.
We refer to \cite[Lemma 3.25]{ChuLas} for a similar assertion.
\end{proof}
Now we complete the proof of Proposition~\ref{pr:q-st}.
Using (\ref{K-damp4-a}) with $\th\in  (0,1]$ we obtain that
\begin{equation}\label{dt-zt}
|(D(t),z_t)|\le C_{R,\eps}|z_t|_\th^2+\eps |z|^2_1
(1+|u^1_t|_\th^2+|u^2_t|_\th^2)
\end{equation}
for any $\eps>0$. We also have from (\ref{K-damp4-a}) for $\th<1$
and from (\ref{K-damp4-b}) for $\th=1$ that
\[
|(D(t),z)|\le C_{R,\eps}|z_t|_\th^2+\eps |z|^2_1
(1+|u^1_t|_\th^2+|u^2_t|_\th^2)+C_{R,\eps}|z|^2.
\]
The subcritical estimate in (\ref{8.1.1c}) yields
\[
|\Psi_T|\le C_R \int_0^T|z_t|_\th^2dt +\eps  \int_0^T|z|_1^2dt+
 C_{R,\eps,T}\int_0^T|z|^2dt~~~\mbox{for every}~~ \eps>0.
\]
Therefore, Lemma~\ref{le:new-dim1} implies
\begin{equation}
 TE_z(T) +
\int_{0}^{T} E_z(t) dt
\le C_{R,\eps} \int_0^T |z_t(t)|_\th^2 dt+  \eps
 \int_0^T  |z|^2_1
(|u^1_t|_\th^2+|u^2_t|_\th^2)dt
 +  C_{R,\eps,T}\int_0^T|z|^2dt
\label{new-dim1a}
\end{equation}
for  every $T\ge T_0$. By Lemma~\ref{le:max}
using  (\ref{.u-int-bound}) and (\ref{dt-zt})   we have that
\begin{equation*}
    \underset{[0,T]}{\max}E_z(t) \leq 2c_0 E_z(T)
  + C_R\left[  \int_0^T|\cA^{\th/2}z_t|^2dt +   \int_0^T\big| \cA^{1/2}z \big|^2dt\right],
\end{equation*}
where the constants $c_0$ and $C_R$ are independent of $T$.
Therefore, from (\ref{new-dim1a})
\[
TE_z(T) +
\int_{0}^{T}\!\! E_z(t) dt
\le C_{R,\eps} \int_0^T \!\! |z_t(t)|_\th^2 dt+  \eps\
 \underset{[0,T]}{\max}E_z(t)
 \int_0^T\!\!
(|u^1_t|_\th^2+|u^2_t|_\th^2)dt
 +  C_{R,\eps,T}\int_0^T\!\!|z|^2dt
\]
which, after an appropriate choice of $\eps$, implies that
\begin{equation}
 TE_z(T) +
\int_{0}^{T} E_z(t) dt
\le C_{R} \int_0^T |z_t(t)|_\th^2 dt +  C_{R,T}\int_0^T|z|^2dt
\label{new-dim1b}
\end{equation}
for  every $T\ge T_0$. Using (\ref{K-damp1a}), (\ref{8.1.1c}) and (\ref{abs-dif-z_t-0})
we conclude that there exists $\tilde \ga_R>0$ such that
\begin{eqnarray*}
\tilde\ga_R \int_0^T |z_t(s)|_\th^2 dt &\le & E_z(0)-E_z(T)+
 \eps \int_0^T |z(t)|_1^2 dt   \\ &&
+C_{R,\eps} \int_0^T  |z|^2_1
(|u^1_t|_\th^2+|u^2_t|_\th^2)dt
 +  C_{R,\eps}\int_0^T|z|^2dt.
\end{eqnarray*}
Consequently, choosing $\eps$ small enough,  by (\ref{new-dim1b}) we have that
\begin{eqnarray*}
 TE_z(T)  &\le & c_R[ E_z(0)-E_z(T)]
\\ && + \,
C_{R} \int_0^T   E_z(t)
(|u^1_t|_\th^2+|u^2_t|_\th^2)dt
 +  C_{R,T}\int_0^T|z|^2dt.
\end{eqnarray*}
Thus
\[
 E_z(T)
\le \kappa_R E_z(0) +
C_{R} \int_0^T   E_z(t)
(|u^1_t|_\th^2+|u^2_t|_\th^2)dt
 +  C_{R,T}\int_0^T|z|^2dt
\]
with $\kappa_R<1$ and $T\ge T_0$.
Now the standard argument (cf., e.g., \cite[p. 62]{ChuLas} or
\cite[p.414]{cl-book}) leads to (\ref{q-stab-est}).
This concludes the proof of Proposition~\ref{pr:q-st}.
\subsection{Completion of the proof of Theorem~\ref{th:attractor}}
{\bf 1.} Proposition~\ref{pr:q-st} means that the system $(\cH,S(t))$
is quasi-stable in the sense of Definition 7.9.2~\cite{cl-book}.
Therefore, by Proposition~7.9.4~\cite{cl-book}
 $(\cH,S(t))$ is asymptotically smooth. 
By Proposition~\ref{pr:generation} $(\cH,S(t))$ is a gradient system.
Thus, Remark~\ref{re:th2.2}
and Theorem~\ref{Theorem 2.2.} imply that there exists a compact global
attractor. By the standard results on gradient systems with compact attractors
 (see, e.g., \cite{BabinVishik,Chueshov,Temam}) we have that $\Ac=\mathbb{M}_+(\cN)$
and (\ref{n-stab}) holds.
\smallskip\par\noindent
{\bf 2.} Since  $(\cH,S(t))$
is quasi-stable,   the finiteness of the fractal dimension ${\rm dim}_f\Ac$
follows from  Theorem~7.9.6~\cite{cl-book}.
\smallskip\par\noindent
{\bf 3.}  To obtain the result on regularity stated in (\ref{u-smth}) and
(\ref{u-smth2}), we apply  Theorem~7.9.8~\cite{cl-book}.
\smallskip\par\noindent
{\bf 4.}  One can see from (\ref{abs-1}) and Theorem~\ref{pr:wp1}
that any weak solution $u(t)$ possesses the property
\[
\int_t^{t+1}|\:u(\tau)|^2_{-1}d\t \le C_R~~\mbox{for all}~~ t>0,
\]
provided $(u_0;u_1)\in B_R=\{y\in\cH\, : \; |y|_\cH\le R\}$.
This implies that $t\mapsto S(t)y$ is a $1/2$-H\"older continuous function
with values in $\tilde\cH =H_\th\times H_{-1}$ for every $y\in B_R$.
Therefore, the existence of a fractal exponential attractor follows from
 Theorem~7.9.9~\cite{cl-book}.
\smallskip\par\noindent
{\bf 5.}
To prove the statement concerning determining functionals,
we use the same idea as in the proof of
 Theorem~8.9.3~\cite{cl-book}, see also  Theorem~7.9.11~\cite{cl-book}.

\subsection{Proof of Theorem~\ref{th:rate}}
We use the same idea as  in \cite{ChuLas_JDDE_2004,ChuLas,cl-book}.
\par 
Since $\cN$ is finite by (\ref{n-stab}) in Theorem~\ref{th:attractor} we have that
for any $y\in \cH$ there exists
  an equilibrium $e=(\phi;0)\in\cN$ such that
\begin{equation}\label{8.6.3}
| S(t)y -e|_\cH\to 0,~~ t\to\infty.
\end{equation}
Thus we need only to prove that $S(t)y$ tends to $e$ with the stated rate.
\par
Let $S(t)y=(u(t);u_t(t))$.
We can assume that $\sup_{t\ge 0} | S(t)y|_\cH\le R$,
for some $R>0$.
The function $z(t)= u(t)- \phi$ satisfies the following
equation
\begin{equation}\label{8.6.4}
z_{tt}(t)  + D(\phi+z(t), z_t(t) )+{\cal A} z(t)+ F(\phi +z(t))-F(\phi)=0,\; t>0.
\end{equation}
Let 
$\widetilde{\cE}(t) = E_z(t)+ \Phi(t)$,
where
$E_z(t)$ is the same as in Lemma~\ref{le:new-dim1} and 
$$
\Phi (t)=  \Pi(\phi+z(t))-\Pi (\phi)-(F(\phi), z)
\equiv \int^1_0 (F(\phi +\lambda z)-F(\phi),z)d\lambda.
$$
One can see that
\begin{equation}\label{8.6.7}
\widetilde{\cE}(t) +  \int_0^t
(D(\phi+z(\t),z_t(\tau)),z_t(\tau))d\tau= \widetilde{\cE}(0).
\end{equation}
In particular, we have that $\widetilde{\cE}(t)$ is non-increasing.
Moreover, since $(z;z_t)\to 0$ in $H_1\times H$ as $t\to+\infty$, we
have that $\widetilde{\cE}(t)\to 0$ when  $t\to+\infty$.
Thus  $\widetilde{\cE}(t)\ge 0$ for all $t\ge 0$.
It is also clear from (\ref{8.1.1c}) that
\begin{equation}\label{8.6.8}
|\widetilde{\cE}(t) - E_z(t)| \le C_R |z(t)|_{1-\delta} |z(t)|_\th
\le \eps \vert z(t) \vert_1^{2}+ C_{R,\eps} |z(t)|^2,\quad \forall\, \eps>0.
\end{equation}
\par
Applying (\ref{new-dim1b}) for  the case when $u^1(t)=u(t)$ and $u^2(t)=\phi$
we obtain that
\begin{equation}
 TE_z(T) +
\int_{0}^{T} E_z(t) dt
\le C_{R} \int_0^T |z_t(t)|_\th^2 dt +  C_{R,T}\max_{[0,T]}|z(t)|^2
\label{new-dim1b-st}
\end{equation}
for $T\ge T_0$ with some $T_0>0$.
Now we prove the following lemma.

\begin{lemma}\label{le8.6.1}
Let $z(t)$ be a weak solution to (\ref{8.6.4}) such that
\begin{equation}\label{8.6.10}
\int_{T-1}^T E_z(t)dt\le\delta~~\mbox{and}~~
\sup_{t\in \R_+}E_z(t)\le \varrho
\end{equation}
with some $\delta,\varrho>$ and $T>1$.
Then there exists $\delta_0>0$ such that
\begin{equation}\label{8.6.11}
\max_{[0,T]}|z(t)|^2\le
C   \int_0^T |z_t(t)|_\th^2 dt
\end{equation}
for every $0<\delta\le\delta_0$,
where the constant $C$ may depend on $\delta$, $\varrho$ and $T$.
\end{lemma}
\begin{proof} Assume that (\ref{8.6.11}) is not true. Then for some
$\delta>0$ small enough
  there exists a sequence
of solutions  $\{z^n(t)\}$ satisfying (\ref{8.6.10}) and
such that
\begin{equation}\label{8.6.11a}
\lim_{n\rightarrow \infty}\left\{\max_{[0,T]}|z(t)|^2
\left[
 \int_0^T |z^n_t(t)|_\th^2 dt\right]^{-1}\right\}
=\infty.
\end{equation}
By (\ref{8.6.10}) $\max_{[0,T]}|z^n(t)|^2\le C_\varrho$ for all $n$.
Thus (\ref{8.6.11a}) implies that
\begin{equation}\label{8.6.12}
\lim_{n\rightarrow \infty} \int_0^T |z^n_t(t)|_\th^2 dt
=0.
\end{equation}
Therefore we can assume that there exists $z^*\in H_1$
such that
\begin{equation}\label{8.6.14}
(z^n; z^n_t)\to (z^*;0)~~\mbox{ $*$-weakly in}~~
L_\infty(0,T; H_1\times \cH).
\end{equation}
It follows
from (\ref{K-damp1}) and (\ref{8.6.12})   that
for $u^n(t)=\phi+z^n(t)$ we have the relation 
\[
 \lim_{n\to\infty}  \int_0^T |( D(u^n(t),u^n_t(t)), \psi(t) )| dt=0~~ 
\mbox{ for any $\psi\in L_2(0,T; H_\th)$.}
\]
This allows us to conclude that  $u^*=\phi+z^*\in H_1$ solves
the problem $\cA u+F(u)=0$.
From (\ref{8.6.10}) we have
that $|\cA^{1/2}(u^*-\phi)|^2\le 2\delta$.
If we choose $\delta_0>0$ such that
$|\cA^{1/2}(\phi_1-\phi_2)|^2> 2\delta$ for every couple $\phi_1$ and $\phi_2$
of stationary solutions (we can do it because the set $\cN$ is finite),
then we can conclude that $u^*=\phi$ provided $\delta\le\delta_0$.
Thus we have $z^*=0$ in (\ref{8.6.14}).
\par
Now  we  normalize the sequence $z^n$ by defining
 $\hat{z}^n \equiv c_n^{-1}z^n$ with $c_n=\max_{[0,T]}|z(t)|$,
where we account only for a suitable subsequence of nonzero terms in $c_n$.
It is clear from (\ref{8.6.14}) with $z^*=0$ that $ c_n \rightarrow 0$ as  $n\rightarrow \infty$.
By (\ref{8.6.11a})  we also have that
\begin{equation}\label{hat-z-cnv}
\int_0^T| \hat{z}_t^n(t)|_\th^2dt\to 0~~\mbox{as}~~
n\to\infty.
\end{equation}
Relations (\ref{new-dim1b-st}) and (\ref{hat-z-cnv})
imply the following uniform estimate
\[
\sup_{t\in [0,T]}\left\{| \hat{z}_t^n(t)|^2 +
|\cA^{1/2} \hat{z}^n(t)|^2\right\}\le C,~~
n=1,2,\ldots.
\]
Thus we can suppose that there exists $\hat{z}^*\in H_1$
such that
\begin{equation}\label{8.6.14a}
(\hat{z}^n; \hat{z}^n_t)\to (\hat{z}^*;0)~~\mbox{ $*$-weakly in}~~
L_\infty(0,T; H_1\times H ).
\end{equation}
The function $\hat{z}^n$ satisfies the equation
\begin{equation}\label{8.6.15}
  \hat{z}^n_{tt} + \frac{1}{c_n}D(\phi+z^n, z^n_t)+\cA \hat{z}^n +
\frac{1}{c_n} [F(\phi+z^n_t) - F(\phi)]=0.
\end{equation}
As above,
from (\ref{K-damp1}) and (\ref{hat-z-cnv}) we  conclude   that
\[
\frac{1}{c_n} \int_0^T |( D(\phi+z^n, z^n_t), \psi )| dt\to 0~~
\mbox{as}~~n\to\infty ~~\mbox{for any $\psi\in L_2(0,T; H_\th)$.}
\]
It  also follows from (\ref{8.6.1a}) and (\ref{8.6.1b}) that
\[
\frac{1}{c_n} [F(\phi+z^n) - F(\phi)]\to \langle F'(\phi), \hat{z}^*\rangle
~~\mbox{weakly in}~~
L_2(0,T; H_{-1}).
\]
Therefore,  after the limit transition in (\ref{8.6.15})
we conclude that $\hat{z}^*$ satisfies
$\cA\hat{z}^*+ \langle F'(\phi), \hat{z}^*\rangle=0$ and, by
hyperbolicity of $\phi$ we conclude that
$\hat{z}^* =0$. Thus (\ref{8.6.14a}) and Aubin-Dubinski theorem 
\cite[Corollary 4]{Simon} imply
that
$\max_{[0,T]}|\hat{z}^n|\to 0$   as $n\to\infty$,
which is impossible.
\end{proof}
\medskip\par\noindent
{\bf Completion of the proof of Theorem~\ref{th:rate}.}
By (\ref{8.6.3}) we choose $T_0>0$ such that (\ref{8.6.10})
holds with $\delta\le \delta_0$ and $T>T_0$.
From (\ref{8.6.8}) we have that $\widetilde{\cE}(T)\le C_R E_z(T)$.
The energy relation in (\ref{8.6.7}) and the lower bound in 
(\ref{K-damp1}) yield that
\begin{equation}\label{e-tt}
\int_0^T |z_t(t)|_\th^2 dt\le   C_R \left[\widetilde{\cE}(0) -
\widetilde{\cE}(T)\right].
\end{equation}
Therefore Lemma~\ref{le8.6.1}, relation (\ref{new-dim1b-st})
 and the energy relation in (\ref{8.6.7}) imply that
$\widetilde{\cE}(T)
\leq \gamma_R  \widetilde{\cE}(0)$ for some $0<\gamma_R<1$. 
This implies that  $\widetilde{\cE}(mT)
\leq \gamma_R^m  \widetilde{\cE}(0)$ for $m=1,2,\ldots$
By (\ref{8.6.8}), (\ref{8.6.11}) and (\ref{e-tt}) we have that 
\[
E_z(mT)\le 2\widetilde{\cE}(mT)+ C_R\max_{[mT,(m+1)T]}|\hat{z}^n|^2
\le C_R\widetilde{\cE}(mT),\quad m=1,2,\ldots
 \]
Thus $E_z(mT)\le C_R \gamma_R^m$ for $m=1,2,\ldots$
Now using (\ref{dif-bnd}) we obtain (\ref{8.6.2}). The proof is complete.

\section{Applications}\label{s:appl}
As it was mentioned in the introduction,
our main applications are plate models.

\subsection{Plate models}
For
the definiteness, we concentrate on the hinged boundary conditions
(the results remain true with other types of self-adjoint
boundary conditions).
Below $\|\cdot\|_s$ is the norm in the Sobolev space
$H^s(\O)$ of order $s$.

\smallskip\par
{\sc Forcing term:} We first check that the forcing term $F$
satisfies Assumption~\ref{A1}(F) for all cases described above.
\par
In the case of the Kirchhoff model, the  embeddings
\[
H^s(\O)\subset L_{2/(1-s)}(\O),\quad
 L_{2/(1+s)}(\O)\subset H^{-s}(\O),\quad
H^{1+\eta}(\O)\subset L_{\infty}(\O)
\]
for $0<s<1$ and $\eta>0$ imply that for any $\th\in [1/2,1]$ the force
$F$   given by (\ref{kirch-force})   
satisfies (\ref{8.1.1c}) with some $\delta<1/2$. 
In the case $\th\in (0,1/2)$ we rely on the inequality
\[
\|f\cdot g\|_s\le C\| f\|_{s+\s} \|g\|_{1-\s}, \quad f,g\in H^1(\O), 
\]
which holds  for any $0<s<1$ and $0<\s<1-s$ (see \cite[Lemma 1.4.1]{cl-book})
and implies that
\[
\|f_1\cdots f_l\|_s\le C\prod_{i=1}^l \| f_i\|_{1-\delta}, \quad 0<s<1,
\]
where $0<\delta<1$ depends  on $s$ and  $l$.  
Using the latter inequality we can find that
\begin{eqnarray*}
\|{\rm div}\left\{ |\nabla u_1|^m\nabla u_1- |\nabla u_2|^m\nabla u_2\right\}\|_{-2\th} &\le & C
\| |\nabla u_1|^m\nabla u_1- |\nabla u_2|^m\nabla u_2 \|_{1-2\th}
\\
&\le & C\|u_1-u_2\|_{2-\delta}\left( \|u_1\|^m_2+\|u_2|^m_2\right),
\end{eqnarray*}
where $m$ is either $q$ or $r$.
Thus (\ref{8.1.1c}) holds for every $0< \th\le 1$.
We also have
\[
\Pi(u)=\int_\O \Phi(u(x))dx +\frac{\kappa}{q+2}\int_\O |\nabla u(x)|^{q+2}dx
-\frac{\kappa\mu}{r+2}\int_\O |\nabla u(x)|^{r+2}dx 
 -\int_\O u(x) p(x)dx,
\]
where $\Phi(s)=\int_0^s\f(\xi)d\xi$ is the antiderivative of $\f$.
It follows from
\eqref{phi_condition}
that there exist $\ga<\l^2_1$ and $C\ge 0$ such that
$\Phi(s)\ge -\ga s^2/2-C$ for all $s\in\R$.
This implies \eqref{8.1.1c1}.
\par
In the case  of the von Karman model, we have that the Airy
stress function $v(u)$ defined in (\ref{airy-1}) satisfies the inequality
\[
\| [u_1,v(u_1)]-  [u_2,v(u_2)]\|_{-\eta}\le
C(\|u_1\|_2^2+\|u_2\|_2^2)\| u_1-  u_2\|_{2-\eta}
\]
for every $\eta\in [0,1]$ (see Corollary 1.4.5  in \cite{cl-book}).
Thus, (\ref{8.1.1c}) holds for every $0<\th\le 1$ with $\delta=\theta$.
The potential energy $\Pi$ has the form
\[
\Pi(u)=\frac14\int_\O\left[ |v(u)|^2 -2([u,F_0]-2 p) u\right] dx
\]
and possesses the properties listed in Assumption~\ref{A1}(F),
see, e.g., \cite[Chapter 4]{cl-book}.
\par
One can also see that the Berger model satisfies
 Assumption~\ref{A1}(F) for every $0<\th\le 1$,
and (\ref{8.1.1c}) holds  with $\delta=\theta$;
for some details see \cite[Chapter 4]{Chueshov} and \cite[Chapter 7]{ChuLas}.
\smallskip\par
{\sc Damping terms:} Now we consider possible forms of the damping operator.
\par
{\bf Case $\th=1$:}
  Our main example is (\ref{damping}) under the following
conditions:
\begin{itemize}
\item $\s_0(s_1)$, $\s_1(s_1,s_2,s_3)$ and $g(s_1,s_2)$
are locally Lipschitz functions of $s_i\in\R$, $i=1,2,3$;
\item
 $\s_0(s_1)>0$, $\s_1(s_1,s_2,s_3)\ge 0$ and $g(s_1,s_2)s_2\ge 0$
for all $s_i$;
\item  there exists $q\ge 0$ such that
\[
|\s_1(\xi)-\s_1(\xi^*)|\le C_R |\xi-\xi^*| (1+|\xi_2|^q+|\xi_3|^q+|\xi_2^*|^q+|\xi_3^*|^q)
\]
for all $\xi=(\xi_1,\xi_2,\xi_3),\, \xi^*=(\xi^*_1,\xi^*_2,\xi^*_3)\in \R^3$
 such that
$|\xi_1|, |\xi_1^*|\le R$;
\item  there exist $q_1\le 4$ and $q_2\le 2$ such that
\begin{equation}\label{g-gr}
|g(\xi)-g(\xi^*)|\le C_R\left[ (1+|\xi_2|^{q_1}+|\xi_2^*|^{q_1})
|\xi_1-\xi_1^*|
+ (|\xi_2|^{q_2}+|\xi_2^*|^{q_2}) |\xi_2-\xi_2^*|\right]
\end{equation}
for all $\xi=(\xi_1,\xi_2),\, \xi^*=(\xi^*_1,\xi^*_2)\in \R^2$ such that
$|\xi_1|, |\xi_1^*|\le R$.
\end{itemize}
 In particular, we can take
\[
\s_1(u,\nabla u)=\s_{10}(u)+\s_{11}(u)|\nabla u|^r~~\mbox{and}~~
g(u,u_t)=g_0(u)  u_t +g_1(u)u_t^3,
\]
where $\s_{10}$, $\s_{11}$, $g_0$ and $g_1$ are  nonnegative 
locally Lipschitz functions, $r\ge 1$.
\smallskip\par

{\bf Case $\th=1/2$:} We consider the case when $\s_0\equiv 0$ and
$\s_1$ is independent of $\nabla u$, i.e., we consider
  a damping of the form
\begin{equation}\label{damping-ss}
D(u,u_t)=
-{\rm div}\,\left[\s_1(u)\nabla u_t\right]+g(u,u_t).
\end{equation}
One can see  that Assumption~\ref{A1}(D) holds with $\th=1/2$
under the following conditions:
\begin{itemize}
\item $\s_1(s_1)$ and $g(s_1,s_2)$
are locally Lipschitz functions of $s_i\in\R$, $i=1,2$;
\item
 $\s_1(s_1)>0$ and $g(s_1,s_2)s_2\ge 0$
for all $s_i$;
\item  the function $g$ satisfies (\ref{g-gr}) with some $q_1<3$ and $q_2<2$.
\end{itemize}\smallskip\par
We can also consider more general (anisotropic) damping operators
$D(u,u_t)$ which are defined variationally  by the formula
\[
(D(u,u_t),\psi)=\sum_{i,j,k,l}\int_\O a_{ijkl}(u)u_{x_ix_j t}\psi_{x_ix_j} dx
+\sum_{i,j}\int_\O b_{ij}(u)u_{x_it}\psi_{x_j} dx +
\int_\O c(u,\nabla u, u_t) \psi dx
\]
for every $\psi\in (H^2\cap H^1_0)(\O)$ under appropriate positivity  and
smoothness hypotheses concerning the coefficients.

\subsection{Wave equation with strong damping}
As an example, we  can also   consider
the following wave equation on a bounded domain $ \O$ in $\R^3$
 with a nonlocal damping coefficient:
\begin{equation}\label{wave-str}
u_{tt}- \s_0(\|u\|_\eta) \D u_t + \s_1(u)u_t  -\D u+ \f(u)=f(x),\quad u\big|_{\d\O}=0.
\end{equation}
We assume that $\eta<1$ and the following conditions
concerning the damping functions  $\s_0$ and $\s_1$ are valid:
(i) $\s_j(s)$
are locally Lipschitz functions of $s\in\R$, $j=0,1$;
(ii) $\s_0(s)>0$, $\s_1(s)\ge 0$ for all $s\in \R$;
(iii)
there exists $q_1< 3$  such that
\[
|\s_1(\xi)-\s_1(\xi^*)|\le C\left(1+|\xi|^{q_1}+|\xi^*|^{q_1}\right)
|\xi-\xi^*|,\quad \xi,\xi^*\in \R.
\]
We also assume that
the source term $\f\in C^1(\R)$
possesses the properties
\[
  \underset{|s|\to\8}{\liminf}\left\{\f(s)s^{-1}\right\}>-\l_1,\quad
|\f'(s)|\le C(1+|s|^q),\; s\in\R,\quad q<4,
\]
 where $\l_1$ is the first eigenvalue of the
Laplacian with the Dirichlet boundary conditions.
\par
In this case we can apply Theorem~\ref{th:attractor} for the model in (\ref{wave-str}) with $\th=1$.
We also note that basing on a requirement like (\ref{new-req}) we 
can also cover the case when $\eta=1$ in (\ref{wave-str}). As for the
case of the critical growth exponents ($q_1=3$ and $q=4$) of the damping 
coefficient $\s_1$ and the force $\f$, we cannot apply here our abstract
 approach. This case requires a separate consideration involving
a  specific structure of the model.
\par 
In a similar way, we can also consider the wave model   (\ref{wave-str}) in
arbitrary dimension $d$ and with another structure of the damping operator.

\end{document}